\documentclass{IET}

\usepackage{multirow}
\usepackage{textcomp}
\usepackage{xcolor}
\usepackage{fancyhdr}
\makeatletter

 \let\oldforeign@language\foreign@language
 \DeclareRobustCommand{\foreign@language}[1]{%
   \lowercase{\oldforeign@language{#1}}}
\newenvironment{lyxlist}[1]
{\begin{list}{}
{\settowidth{\labelwidth}{#1}
 \setlength{\leftmargin}{\labelwidth}
 \addtolength{\leftmargin}{\labelsep}
 }}
{\end{list}}

\newcommand{\splitatcommas}[1]{%
	\begingroup
	\begingroup\lccode`~=`, \lowercase{\endgroup
		\edef~{\mathchar\the\mathcode`, \penalty0 \noexpand\hspace{0pt plus 0em}}%
	}\mathcode`,="8000 #1%
	\endgroup
}

\renewcommand{\fnum@figure}{Fig.~\thefigure}

\makeatother

\begin{document}

\title{Optimal Allocation of Series FACTS Devices in Large Scale Systems\newline\vskip3pt}

\author[1,*]{Xiaohu Zhang}
\affil{Department of Electrical Engineering and Computer Science, The University of Tennessee, Knoxville, TN 37996, USA}

\author[1]{Kevin Tomsovic}
%\author[1]{Bailu Xiao}

\author{Aleksandar Dimitrovski}
\affil{Department of Electrical and Computer Engineering, University of Central Florida, FL 32816, USA}
\affil[*]{xzhang46@vols.utk.edu}

\abstract{Series FACTS devices, such as the variable series reactor (VSR), have the ability to continuously regulate the transmission line reactance so as to control power flow. This paper presents a new approach to optimally locating such devices in the transmission network considering multiple operating states and contingencies. To investigate optimal investment, a single target year planning with three different load patterns is considered. The transmission contingencies may occur under any of the three load conditions and the coupling constraints between base case and contingencies are included. A reformulation technique transforms the original mixed integer nonlinear programming (MINLP) model into mixed integer linear programing (MILP) model. To further relieve the computational burden and enable the planning model to be directly applied to practical large scale systems, a two phase decomposition algorithm is introduced. Detailed numerical simulation results on IEEE 118-bus system and the Polish 2383-bus system illustrate the efficient performance of the proposed algorithm.}

\maketitle

\section*{Nomenclature}
\subsection*{Indices}
\begin{lyxlist}{00.00.0000}
\item[$i, \ j$] Index of buses.
\item[$k$] Index of transmission elements.
\item[$n$] Index of generators.
\item[$m$] Index of loads.
\item[$c$] Index of states; $c=0$ indicates the base case; $c>0$ is a contingency state.
\item[$t$] Index of load levels.
\end{lyxlist}

\subsection*{Variables}
\begin{lyxlist}{00.00.0000}
\item[$P^g_{nct}$] Active power generation of generator $n$ for state $c$ at load level $t$.
\item[$P_{kct}$] Active power flow on branch $k$ for state $c$ at load level $t$.
\item[$x^V_{k}$] Reactance of a VSR at branch $k$.
\item[$\theta_{kct}$] The angle difference across the branch $k$ for state $c$ at load level $t$.
\item[$\delta_{k}$] Binary variable associated with installing a VSR on branch $k$.
\item[$\Delta P^{g,up}_{nct}$] Active power generation adjustment up of generator $n$ for state $c$ at load level $t$.
\item[$\Delta P^{g,dn}_{nct}$] Active power generation adjustment down of generator $n$ for state $c$ at load level $t$.
\item[$\Delta P^d_{mct}$] Load shedding quantity of load $m$ for state $c$ at load level $t$.
\end{lyxlist}

\subsection*{Parameters}
\begin{lyxlist}{00.00.0000}
%\item[$b_k$] Susceptance for branch $k$.
\item[$a_n^g$] Cost coefficient for generator $n$.
\item[$a_n^{g,up}$] Cost coefficient for generator $n$ to increase active power.
\item[$a_n^{g,dn}$] Cost coefficient for generator $n$ to decrease active power.
\item[$a_{LS}$] Penalty for the load shedding.
\item[$A_h$] Annual operating hours: 8760 h.
\item[$A_I$] Annual investment cost for VSR.
\item[$N_{kct}$] Binary parameter associated with the status of branch $k$ for state $c$ at load level $t$. 
\item[$\theta_k^{\max}$] Maximum angle difference across branch $k$: $\pi/3$ radians.
\item[$P_{nct}^{g,\min}$] Minimum active power output of generator $n$ for state $c$ at load level $t$.
\item[$P_{nct}^{g,\max}$] Maximum active power output of generator $n$ for state $c$ at load level $t$.
\item[$P_{mct}^d$]  Active power consumption of demand $m$ for state $c$ at load level $t$.
\item[$S_{kct}^{\max}$] Thermal limit of branch  $k$ for state $c$ at load level $t$. 
%\item[$\theta_k^{\min}$] Minimum angle difference across branch $k$.
\item[$R_n^{g,up}$] Ramp up limit of generator $n$ .
\item[$R_n^{g,dn}$] Ramp down limit of generator $n$.
\item[$\pi_{ct}$] Duration of state $c$ at load level $t$.
\end{lyxlist}
\subsection*{Sets}
\begin{lyxlist}{00.00.0000}
\item[$\mathcal{B}$] Set of buses. 
\item[$\mathcal{B}_{ref}$] Set of reference bus.
\item[$\Omega_{L}$] Set of transmission lines.
\item[$\Omega_{L}^i$] Set of transmission lines connected to bus $i$.
\item[$\mathcal{G}$] Set of on-line generators.
\item[$\mathcal{G}_{i}$] Set of on-line generators located at bus $i$.
\item[$\mathcal{G}_{re}$] Set of on-line generators that allow to rescheduling.	
\item[$\mathcal{D}$] Set of loads.
\item[$\mathcal{D}_i$] Set of loads located at bus $i$.
\item[$\Omega_{0}$] Set of base operating states.
\item[$\Omega_{c}$] Set of contingency operating states.
\item[$\Omega_{V}$] Set of candidate transmission lines to install VSR.
\item[$\Omega_{T}$] Set of load levels.
\end{lyxlist}
Other symbols are defined as required in the text.
\section{Introduction}\label{sec1}

Due to the power market restructuring and the rapid introduction of renewables, the electric power industry is going through profound changes across technical, economic and organizational concerns. While deregulation has been able to deliver on some promises, such as, the reduction in electricity prices and new innovative technologies to improve the grid efficiency, it has also led to strains on the transmission system, which was not designed for this new structure. Increasing electricity consumption, less predictable power flows and extensive adoption of renewable energy have led to increasing power grid congestion \cite{mybibb:ps_economics}. 
One option to relieve congestion in the transmission network is through power system expansion, which involves building new power plants and transmission lines in critical areas. This option suffers from the difficulties in obtaining right-of-way, high cost and long construction times. An alternative to increase the effective transmission capability is by installing the power flow control devices such as the Flexible AC Transmission Systems (FACTS) in the existing network \cite{mybibb:FACTS1,mybibb:optfacts,mybibb:FACTS2}.

A series FACTS controller, variable series reactor (VSR), has the ability to vary the effective transmission line reactance so it is suitable for power flow control. Increasing the impedance on the congested lines can shift power to underutilized transmission lines nearby while decreasing the impedance of the transmission line can increase power transferred on that line assuming thermal limits have not been reached. With the rapid developments in the power electronics technology, VSR offers excellent control and flexibility. Moreover, according to the Green Electricity Network Integration (GENI) program \cite{mybibb:geni}, it is anticipated that there will be new FACTS-like devices with far cheaper cost \cite{mybibb:Aleks} available to be installed in the transmission network across US in the future. Given these considerations, efficient planning models and algorithms, which directly work for practical large scale systems, should be developed to provide an optimal planning for the application of FACTS devices.   

The allocation and utilization of FACTS devices has been studied extensively during the last several decades. \textcolor{black}{Considering the nonlinear and non-convex characteristics of the power flow equations, different heuristic approaches such as differential evolution (DE) \cite{mybibb:de_1}, genetic algorithm (GA) \cite{mybibb:ga_1,mybibb:GA_FACT_market,mybibb:GAFACTS,mybibb:ga_tool_facts}, particle swarm optimization (PSO) \cite{mybibb:PSO_FACTS,mybibb:security_pso} have been leveraged to optimally allocate FACTS devices}. Priority indices are another class of methods for locating FACTS devices. To decide the optimal locations of thyristor controlled series compensator (TCSC), the authors in \cite{mybibb:TCSCsens} define the priority indices as the weighted sensitivities of the system transfer capability with respect to each line reactance. In \cite{mybibb:facts_priority_index}, the difference of locational marginal price (LMP) across a transmission line is computed with a standard optimal power flow (OPF) problem. The locations of TCSC are determined based on the magnitudes of the LMP difference. A main drawback of the priority indices methods is that the quality of the solution regarding optimality cannot be guaranteed.  

With rapid advances in mathematical programming algorithms, these methods have been widely used to analyze the impacts of FACTS devices. The authors in \cite{mybibb:TCSCLFB} leverage line flow based equations proposed by \cite{mybibb:LFB} to determine the locations and compensation levels of TCSC via mixed integer linear programming (MILP) and mixed integer quadratic programming (MIQP). To linearize the product of two continuous variables in the constraints, one variable in the product is relaxed by its upper and lower bound. In addition, only one load pattern is considered so the planning model is suitable for preliminary analysis. In \cite{mybibb:ts_tcsc_wind}, the branch and price algorithm is utilized to co-optimize the locations of transmission switch and TCSC considering the wind power uncertainties. In \cite{mybibb:bd_sopf1,mybibb:bd_sopf2}, Benders Decomposition is used to investigate the benefits of VSR devices in the security constrained optimal power flow (SCOPF) problem. The master problem is to minimize the generation cost with the pre-located FACTS devices and the subproblem is used to check the feasibility for each contingency. To include VSR in the economic dispatch (ED) problem,  researchers have reformulated the nonlinear programming (NLP) model into an MILP model \cite{mybibb:asu_FACTS1,mybibb:Tao_Ding}.  

Reference \cite{mybibb:opf_redispatch_facts} demonstrates that with appropriately installed TCSC and static VAR compensator (SVC) in the power network, the operation cost during contingencies can be reduced. Hence, from the planning point of view, including the contingency constraints in FACTS placement model could provide a more accurate and useful investment strategy for the system planners. However, the task is not trivial since adding contingency constraints increases the model size significantly and leads to excessive computational burdens. Reference \cite{mybibb:pso_sqp_facts_wind} and \cite{mybibb:PSO_SQP_FACTS} adopt the two level hybrid PSO/SQP algorithm to address this problem but the computation burden is large for a small or medium scale system.

This paper proposes a new solution approach to optimally allocate VSR in large scale transmission networks considering the base case and a series of $N-1$ transmission contingencies. We consider a single target year for the planning. Three distinct load levels which denote peak, normal and low load conditions are selected. The original planning model is a large scale mixed integer nonlinear programming (MINLP) model which is difficult to solve for practical systems. A reformulation technique is used to transform the MINLP model into an MILP model. To further relieve the computational burden, a two phase Benders Decomposition separates the problem into base case master problem and a series of subproblems for contingencies. Considering the extensive literature in this area, the contributions of this paper are twofold:
\begin{enumerate}
	\item \textcolor{black}{to develop a planning model to allocate VSR in the transmission network considering a multi-scenario framework and solve the model using mathematical programming rather than the heuristic or sensitivity methods so that the quality of the solution can be ensured;}    
	\item \textcolor{black}{to implement a two phase Benders decomposition for the planning problem which shows high performance even for a practical large scale network considering hundreds of operating conditions.}
\end{enumerate}       

The remaining sections are organized as follows. Section \ref{static_model} illustrates the reformulation technique. The detailed optimization model is given in section \ref{Optimization_model}. In section \ref{solution_BD}, the solution procedure based on Benders Decomposition is described. The IEEE 118-bus system and Polish 2383-bus system are selected for case studies in section \ref{case_study}. Finally, conclusions are given in section \ref{conclusions}.

\section{Reformulation Technique}
\label{static_model}
We leverage the reformulation technique proposed in \cite{mybibb:tep_cvsr,mybibb:investment_naps} to linearize the nonlinear power flow equations due to the introduction of VSR. The procedures are illustrated in this section for completeness. 

\begin{figure}[!htb]
	%\captionsetup{font=footnotesize}
	\centering
	\includegraphics[width=0.6\textwidth]{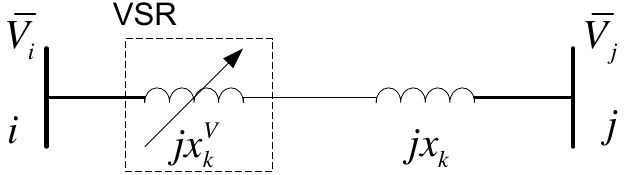}
	\caption{Steady state model of VSR in DCPF.}
	\label{tcsc_steady}
\end{figure}

The static model of a VSR in DC power flow is depicted in Fig. \ref{tcsc_steady}. It can be denoted by a variable reactance $x_k^V$ in series with the line reactance $x_k$.

The resulting susceptance of the transmission line $\tilde{b}_k$ can be expressed as:
\begin{equation}
\tilde{b}_k=-\frac{1}{x_k+x_k^{V}}=-(b_{k}+b_k^{V}) \label{result_susceptance}
\end{equation}
In (\ref{result_susceptance}), $b_k$ can be interpreted as the negative susceptance of line $k$ without VSR and $b_k^V$ is the susceptance change introduced by the VSR. Their expressions are give as: 
\begin{align}
&b_{k}=\frac{1}{x_k}   \\
&b_k^{V}=-\frac{x_k^{V}}{x_k(x_k+x_k^{V})}
\end{align}

The active power transferred on transmission line $k$ is:
\begin{align}
&P_k=(b_{k}+\delta_kb_k^{V})\theta_k  \label{power_CVSR}  \\
&b_{k,V}^{\min} \le b_k^{V} \le b_{k,V}^{\max}   \label{bcvsr_range} 
\end{align}
The binary variable $\delta_k$ is introduced to indicate the installation of a VSR on the transmission line. 

In (\ref{power_CVSR}), the nonlinear term is $\delta_kb_k^{V}\theta_k$. To linearize the nonlinear term, a dummy variable $\psi_k$ is first introduced:
\begin{equation}
\psi_k=\delta_kb_k^{V}\theta_k  \label{wij1}
\end{equation}
Then the power flow equation (\ref{power_CVSR}) can be rewritten as:
\begin{equation}
P_k=b_{k}\theta_k+\psi_k  \label{PCVSR_wij1}
\end{equation}
and we multiply $\delta_k$ with each side of the constraint (\ref{bcvsr_range}) to obtain:
\begin{equation}
\delta_k b_{k,V}^{\min} \le \frac{\psi_k}{\theta_k}=\delta_kb_k^{V} \le \delta_k b_{k,V}^{\max}  \label{if_ineq}
\end{equation}
The allowable range of $\psi_k$ is dependent on the sign of the phase angle difference across transmission line $k$, i.e., $\theta_k$. We introduce a binary variable $y_k$ and use the big-M complementary constraints to linearize (\ref{if_ineq}): 
\begin{align}
&-M_ky_k+\delta_k\theta_kb_{k,V}^{\min} \le \psi_k 
\le \delta_k\theta_kb_{k,V}^{\max}+M_ky_k \label{if_1} \\
%\end{equation}
%\begin{align}
&-M_k(1-y_k)+\delta_k\theta_kb_{k,V}^{\max} \le \psi_k 
\le \delta_k\theta_kb_{k,V}^{\min}+M_k(1-y_k) \label{if_2} 
\end{align}

The physical interpretation of $y_k$ is the flow direction of transmission element $k$. Note that only one of the two constraints (\ref{if_1}) and (\ref{if_2}) will be active during the optimization process and the other one will be always satisfied because of the large number $M_k$. 
%Note that numerical problems occur when $M_k$ is chosen to be too large \cite{mybibb:Tao_bigM,mybibb:TEP_huizhang1}. 

In constraints (\ref{if_1}) and (\ref{if_2}), there still exists a bilinear term $\delta_k\theta_k$ which is the product of a continuous variable and a binary variable. Another variable $v_k=\delta_k\theta_k$ is introduced and linearized by using the big-M method \cite{mybibb:WP_Adams}:
\begin{align}
&-\delta_k\theta_k^{\max} \le v_k \le \delta_k\theta_k^{\max} \label{z1} \\
&\theta_k-(1-\delta_k)\theta_k^{\max} \le v_k \le \theta_k+(1-\delta_k)\theta_k^{\max} \label{z2}
\end{align}
Finally, the constraints (\ref{if_1}) and (\ref{if_2}) can be rewritten as (\ref{if3}) and (\ref{if4}) by replacing $\delta_k\theta_k$ with $v_k$:
\begin{align}
&-M_ky_k+v_kb_{k,V}^{\min} \le \psi_k 
\le v_kb_{k,V}^{\max}+M_ky_k \label{if3} \\
%\end{align}
%\begin{align}
&-M_k(1-y_k)+v_kb_{k,V}^{\max} \le \psi_k\le v_kb_{k,V}^{\min}+M_k(1-y_k)\label{if4} 
\end{align}

Therefore, the nonlinear power flow equations (\ref{power_CVSR}) and (\ref{bcvsr_range}) are reformulated into an MILP format including (\ref{PCVSR_wij1}), (\ref{z1})-(\ref{if4}).

\section{Problem Formulation}
\label{Optimization_model}
With the reformulation, the complete optimization model can be represented as a large scale MILP.
\subsection{Objective Function}
The objective of the proposed planning model is to minimize the total cost for the single target year, which can be formulated as:
\begin{equation}
\min_{\Xi_{\text{OM}}} \ \ \sum_{k\in \Omega_V}A_I\delta_k+\sum_{t\in \Omega_T}(\pi_{0t}C_{0t}+\sum_{c\in \Omega_c}\pi_{ct}C_{ct})  \label{objective}
\end{equation}
There are three components in the objective function (\ref{objective}). Specifically, the first term denotes the annualized investment cost in VSR; the second and third term represent the operation cost under base and contingency states respectively. Under load level $t$, the operation cost for the base state, i.e., $C_{0t}$, can be further expressed as:
\begin{equation}
C_{0t}=\sum_{n\in \mathcal{G}}a_{n}^gP^g_{n0t}  \label{obj_base}
\end{equation}

The operation cost under the contingency state $c$ and load level $t$ comprises four terms:
\begin{align}
C_{ct}=&\sum_{n\in \mathcal{G}}a_{n}^gP^g_{nct}+\sum_{m\in \mathcal{D}}a_{LS}\Delta P^d_{mct}    \nonumber \\
&+\sum_{n\in \mathcal{G}}(a_n^{g,up}\Delta P^{g,up}_{nct}+a_n^{g,dn}\Delta P^{g,dn}_{nct})                        \label{obj_cont}     
\end{align}
The generation cost under each contingency is indicated by the first term; the load shedding cost is denoted by the second term; the generator up and down adjustment costs are represented by the third and fourth term, respectively. Each operating state is associated with a duration time $\pi_{ct}$ so the cost for the individual state is appropriately weighted. Note that the number of operation hours for a target year is 8760, which is give by (\ref{operating_hour}):
\begin{equation}
\sum_{t\in \Omega_T}\pi_{0t}+\sum_{t\in \Omega_T}\sum_{c\in \Omega_c}\pi_{ct}=A_h    \label{operating_hour}
\end{equation}

\subsection{Constraints}
The complete set of constraints are given below from (\ref{norm_c1}) to (\ref{load_shedding}).
\begin{align}
%\big \{
& P_{kct}=N_{kct}b_{k}\theta_{kct},\ k\in \Omega_{L}\backslash \Omega_{V} \label{norm_c1} \\
%& P_{kc}-b_{k}\theta_{kct}-M_{k}'(1-N_{kct}) \le 0,\ k\in \Omega_{L}\backslash \Omega_{V} \label{norm_c2} \\
& P_{kct}=N_{kct}(b_{k}\theta_{kct}+\psi_{kct}),\ k\in \Omega_{V} \label{CVSR_c1}  \\
%& P_{kc}-b_{k}\theta_{kct}-\psi_{kct}-M_{k}'(1-N_{kct}) \le 0,\ k\in \Omega_{V} \label{CVSR_c2}   \\
&-M_{k}y_{kct}+v_{kct}b_{k,V}^{\min} \le \psi_{kct}\le v_{kct}b_{k,V}^{\max}+M_{k}y_{kct},\  k\in \Omega_{V} \label{big_M1} \\
&-M_{k}(1-y_{kct})+v_{kct}b_{k,V}^{\max} \le \psi_{kct}\le  \nonumber \\ & \ \ \ \ \ \ \ \ \ \ \ \ \ \ \ \ \ \ \ \ \ \ v_{kct}b_{k,V}^{\min}+M_{k}(1-y_{kct}),\ k\in \Omega_{V} \label{big_M2}   \\ 
& -\delta_{k}\theta_{k}^{\max} \le v_{kct} \le \delta_{k}\theta_{k}^{\max},\ k \in \Omega_{V} \label{bilinear_c1} \\
& \theta_{kct}-(1-\delta_{k})\theta_{k}^{\max} \le v_{kct}\le \theta_{kct}+(1-\delta_{k})\theta_{k}^{\max},\ k \in \Omega_{V} \label{bilinear_c2} \\
&\sum_{n \in \mathcal{G}_i}P^g_{nct}-\sum_{m \in \mathcal{D}_i}(P^d_{mct}-\Delta P^d_{mct})=\sum_{k\in \Omega_{L}^i } P_{kct} \label{p_bal}  \\ 
& -S_{kct}^{\max} \le P_{kct} \le S_{kct}^{\max}, \ k\in \Omega_L \label{Slimit_E} \\
& P_{nct}^{g,\min}\le P_{nct}^g \le P_{nct}^{g,\max} \label{Pg_limit1} \\
& \theta_{i}=0, \ i\in \mathcal{B}_{ref}     \label{reference_angle}    \\
%&\big \} 
%\end{align}
%\end{align}
%\begin{align}
%\begin{align}
&P_{nct}^g=P_{n0t}^g+\Delta P_{nct}^{g,up}-\Delta P_{nct}^{g,dn}, \ n\in \mathcal{G}_{re}    \label{ramp_gen} \\
&0 \le \Delta P_{nct}^{g,up} \le R_n^{g,up}, \ n\in \mathcal{G}_{re}   \label{ramp_up}   \\
&0 \le \Delta P_{nct}^{g,dn}  \le R_n^{g,dn},\ n\in \mathcal{G}_{re}  \label{ramp_dn} \\
&P_{nct}^g=P_{n0t}^g, \ n\in \mathcal{G}\backslash\mathcal{G}_{re}  \label{g_fix}  \\
&0 \le \Delta P^d_{mct} \le P^d_{mct} \label{load_shedding}
\end{align}
Constraints (\ref{norm_c1})-(\ref{bilinear_c2}), (\ref{Slimit_E})-(\ref{reference_angle}) hold  $\forall c \in \Omega_c\cup \Omega_0, t \in \Omega_T, n \in \mathcal{G}$, constraints (\ref{p_bal}) hold $\forall c \in \Omega_c\cup \Omega_0, t \in \Omega_T, i\in \mathcal{B}$, and constraints (\ref{ramp_gen})-(\ref{load_shedding}) hold $\forall c \in \Omega_c, t \in \Omega_T, m\in \mathcal{D}$.

Constraints (\ref{norm_c1})-(\ref{reference_angle}) are the operating constraints, including base case and contingencies. Specifically, constraint (\ref{norm_c1}) is the power flow on the lines without VSR and constraint (\ref{CVSR_c1}) represents the power flow on the candidate lines to install VSR. We introduce a binary parameter $N_{kct}$ to denote the corresponding status of the transmission element $k$ in state $c$ at load level $t$ \cite{mybibb:TEP_huizhang1}. If $N_{kct}=1$, the line flow equations are forced to hold; otherwise, if the line is in outage, the power flow on that line is forced to be zero. The reformulation considering multiple operating states and load levels are denoted by constraints (\ref{big_M1})-(\ref{bilinear_c2}). Constraints (\ref{p_bal}) ensure the power balance at each bus. The thermal limits of the transmission lines and the active power limits of the generators are considered in (\ref{Slimit_E}) and (\ref{Pg_limit1}). Note that the short term rating for the transmission line is used for the contingency states, which is 10\% higher than the thermal limit under the base operating condition. Finally, constraint (\ref{reference_angle}) sets the bus angle of the reference bus to zero.

Constraints (\ref{ramp_gen})-(\ref{load_shedding}) denote limits under the contingency states. Constraints (\ref{ramp_gen})-(\ref{g_fix}) indicate that only a subset of generators are allowed to redispatch their generation during the contingencies and all the other generators should be fixed at their base operating condition. The load shedding quantity should not exceed the existing load, which is given in  (\ref{load_shedding}).

The optimization variables of the complete planning model from (\ref{objective})-(\ref{load_shedding}) are the elements in set $\splitatcommas{\Xi_{\text{OM}}=\{\theta_{kct},P_{kct},P^g_{nct},\Delta P^d_{mct},\Delta P^{g,up}_{nct},\Delta P^{g,dn}_{nct},\delta_k,y_{kct},v_{kct},\psi_{kct}\}}$.

\section{Solution Approach}
\label{solution_BD}
\textcolor{black}{The size of the MILP model formulated in Section \ref{Optimization_model} dramatically increases with the system size and the number of considered contingencies, which leads to excessive computation. In order to make the optimization model applicable to a practical large system, Benders Decomposition is used to decompose the original optimization model into a master problem and a subproblem. The master problem deals with the base operating condition and the subproblem considers contingencies.} The complicating variables between the master problem and subproblem are $P^g_{n0t}$ and $\delta_k$. 

It should be noted that the prerequisite for Benders Decomposition is that the objective function of the considered problem projected on the subspace of the complicating variables has a convex envelope \cite{mybibb:bd_conejo}. This is not the case in our model due to the existence of the binary flow direction variable $y_{kct}$ in the subproblem. In \cite{mybibb:bingqian_hu}, a modified Benders Decomposition (MBD) is developed for the security constrained unit commitment (SCUC) considering the quick-start units. The main idea is to construct a tighter LP subproblem based on the MILP subproblem and use the tighter LP to generate Benders cuts. We propose an alternative two phase approach in section \ref{solution_procedure}. The simulation results obtained from the proposed approach and MBD are compared in section \ref{case_study}.    

\subsection{Master Problem}
The master problem considers the base operating condition for the three load levels:
\begin{align}
&\min_{\Xi_{\text{MP}}}\ Z_{down}^{(\nu)}=\sum_{t\in \Omega_T}\pi_{0t}C_{0t}^{(\nu)}+\sum_{k\in \Omega_V}A_I\delta_k^{(\nu)}+\alpha^{(\nu)}   \label{master_obj} \\
& \text{subject to:}    \nonumber  \\
&(\ref{norm_c1})-(\ref{reference_angle}) \ \ \text{and}  \nonumber  \\
&\alpha^{(\nu)} \ge \alpha_{down}  \label{accerlerate}  \\
&\alpha^{(\nu)} \ge Z^{(l)}+\sum_{t\in \Omega_T}\sum_{n\in \mathcal{G}}\mu_{nt}^{(l)}(P_{n0t}^{g^{(\nu)}}-P_{n0t}^{g^{(l)}}) \nonumber\\
&\ \ \ \  +\sum_{k\in \Omega_V}\beta_k^{(l)}(\delta_k^{(\nu)}-\delta_k^{(l)}), \ l=1,\cdots,\nu-1  \label{bender_cut}  
\end{align}
Constraints (\ref{master_obj})-(\ref{bender_cut}) hold $\forall c \in \Omega_0, t \in \Omega_T, n\in \mathcal{G}, i\in \mathcal{B}$.

The optimization variables of the master problem are those in the set $\splitatcommas{\Xi_{\text{MP}}=\{\theta_{kct},P_{kct},P^{g}_{nct},\delta_k,y_{kct},v_{kct},\psi_{kct},\alpha\}}$. Note that all the variables are subject to Benders iteration parameter $\nu$. The first and second term in the objective function are the operating cost in the base case and the investment cost for the VSR. $\alpha^{(\nu)}$ denotes the total operating cost during the contingencies. To accelerate the convergence speed, constraint (\ref{accerlerate}) sets a lower bound on $\alpha^{(\nu)}$. Constraint (\ref{bender_cut}) represents the Benders cut, which will be generated once per iteration.

\subsection{Subproblem}
The subproblem for contingency state $c$ and load level $t$ is:
\begin{align}
&\min_{\Xi_{\text{SP}}} Z_{ct}^{(\nu)}=C_{ct}^{(\nu)}+\sum_{i\in \mathcal{B}}h_i(s_{ict,1}^{(\nu)}+s_{ict,2}^{(\nu)})  \label{sub_obj}  \\
&\text{subject to}   \nonumber \\
&(\ref{norm_c1})-(\ref{bilinear_c2}),(\ref{Slimit_E})-(\ref{load_shedding}) \ \ \text{and} \nonumber    \\
%\end{align}
%\begin{align}
&\sum_{n \in \mathcal{G}_i}P^{g^{(\nu)}}_{nct}-\sum_{m \in \mathcal{D}_i}(P^d_{mct}-\Delta P^{d^{(\nu)}}_{mct})  \nonumber \\
&\ \ \ \ \ \ \ \ \ \ \ \ \ \ \ \ \ \ \ +s^{(\nu)}_{ict,1}-s^{(\nu)}_{ict,2}=\sum_{k\in \Omega_{L}^i } P^{(\nu)}_{kct} \label{p_bal1}  \\
&s_{ict,1}^{(\nu)} \ge 0, \ s_{ict,2}^{(\nu)} \ge 0    \label{slack_range}   \\
&P^{g^{(\nu)}}_{n0t}=\hat{P}^g_{n0t} \ \ \ \ : \mu_{nct}^{(\nu)}   \label{fix_P} \\
&\delta_k^{(\nu)}=\hat{\delta}_k \ \ \ \ : \beta_{kct}^{(\nu)}     \label{fix_de} 
\end{align}
Constraints (\ref{sub_obj})-(\ref{fix_de}) hold $\forall c \in \Omega_c, t \in \Omega_T, n\in \mathcal{G}, i\in \mathcal{B}, m\in \mathcal{D}$.

The optimization variables of the subproblem are those in the set $\splitatcommas{\Xi_{\text{SP}}=\{\theta_{kct},P_{kct},P^g_{nct},\Delta P^d_{mct},\Delta P^{g,up}_{nct}, P^{g,dn}_{nct},v_{kct},y_{kct},\delta_k, \psi_{kct},s_{ict,1},s_{ict,2}\}}$. The first term of the objective function is the operating cost in each contingency. Note that although load shedding is allowed in the contingency state, the subproblem can still be infeasible due to generator ramping constraints. Two slack variables $s^{(\nu)}_{ict,1}$ and $s^{(\nu)}_{ict,2}$ are introduced to ensure the subproblem is feasible with $h_i$ a sufficiently large positive constant. The complicating variables are fixed at the value obtained from the master problem in constraints (\ref{fix_P}) and (\ref{fix_de}) while $\mu_{nct}^{(\nu)}$ and $\beta_{kct}^{(\nu)}$ are the dual variables associated with these two constraints.

The sensitivity used to generate Benders cut is the weighted dual variable, which can be expressed as:
\begin{align}
&\mu_{nt}^{(\nu)}=\sum_{c\in \Omega_c}\pi_{ct}\mu_{nct}^{(\nu)}  \label{bd_dual_p}\\
&\beta_k^{(\nu)}=\sum_{c\in \Omega_c}\sum_{t\in \Omega_T}\pi_{ct}\beta_{kct}^{(\nu)}  \label{bd_dual_de}
\end{align}

In the master problem, $Z$ in the Benders cut constraint can be calculated as:
\begin{equation}
Z^{(\nu)}=\sum_{c\in \Omega_c}\sum_{t \in \Omega_T}\pi_{ct}Z_{ct}^{(\nu)}  \label{z_benders_cut}
\end{equation}
With the solution of the subproblem, the upper bound of the objective function for the original problem is calculated as:
\begin{equation}
Z_{up}^{(\nu)}=Z^{(\nu)}+\sum_{t\in \Omega_T}\pi_{0t}\hat{C}_{0t}+\sum_{k\in \Omega_V}A_I\hat{\delta}_k  \label{z_up}
\end{equation}	
The last two terms in (\ref{z_up}) are calculated using the fixed value of $\hat{P}^g_{n0t}$ and $\hat{\delta}_k$.

\subsection{Solution Procedure}
\label{solution_procedure}
\begin{figure}[!htb]
	%\captionsetup{font=footnotesize}
	\centering
	\includegraphics[width=0.55\textwidth]{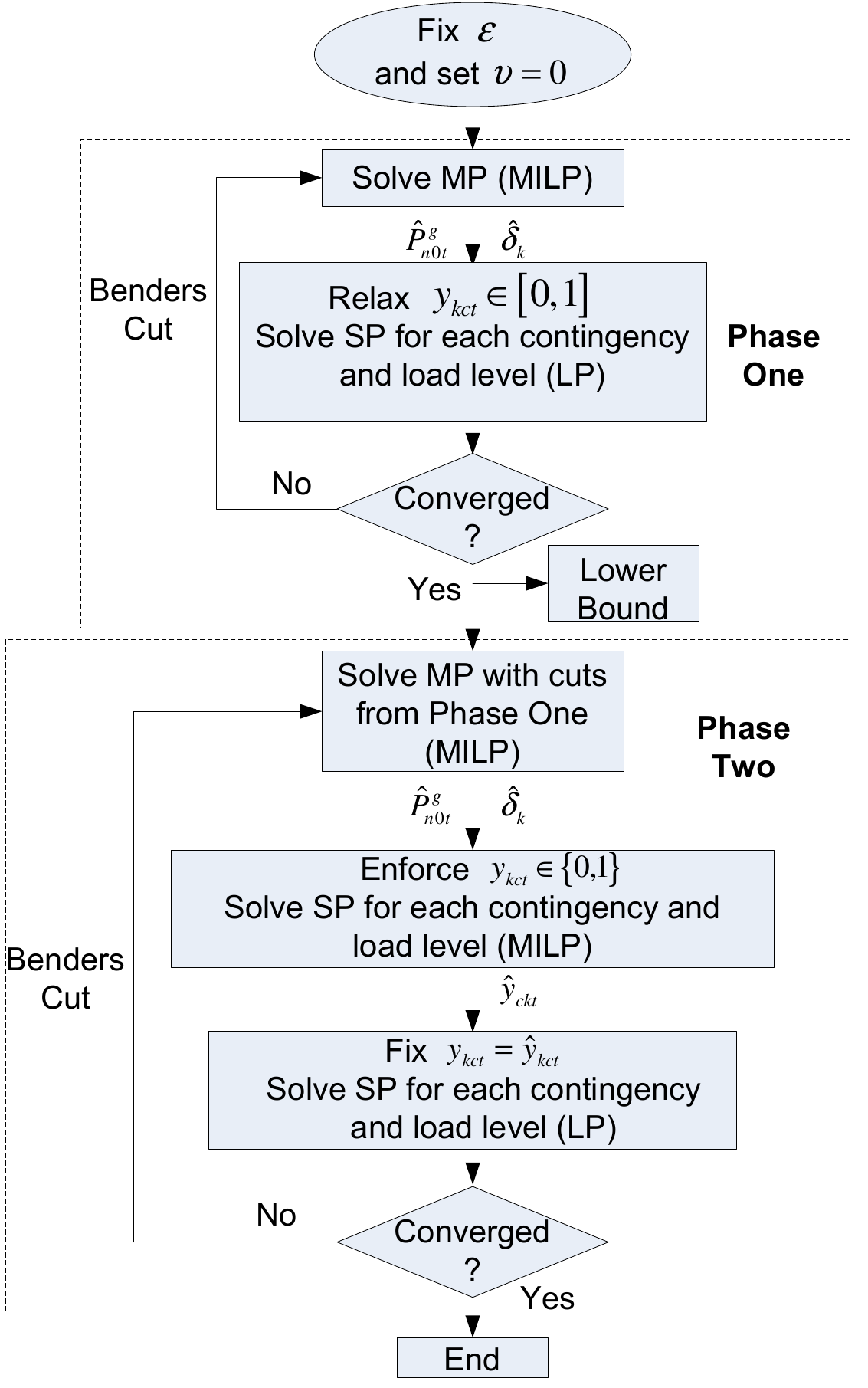}
	\caption{Flowchart of the solution approach.}
	\label{flowchart}
\end{figure}

As mentioned at the beginning of this section, a two phase approach is proposed to solve the planning model. The flowchart of the optimization procedure is shown in Fig. \ref{flowchart}.

The detailed description of the proposed algorithm is given as follows: 
\begin{itemize}
	\item[1):] \textbf{Initialization}: Set a small value $\epsilon$ to control the convergence and initiate the iteration counter $\nu=0$.
	\item[2):] \textbf{Phase one master problem solution}: Solve the master problem considering only the normal operating states. Note that for the first iteration, the master problem is solved without considering any Benders cut, e.g., constraint (\ref{bender_cut}).
	\item[3):] \textbf{Relaxed subproblem solution}: With $P^g_{n0t}$ and $\delta_k$ obtained from the master problem, solve the subproblem as an LP by relaxing the flow direction variable $y_{nct}$ as a continuous variable in $[0,1]$. 
	\item[4):] \textbf{Convergence check}: If $|Z_{up}^{(\nu)}-Z_{down}^{(\nu)}|/|Z_{down}^{(\nu)}|\le \epsilon$, the optimal solution for the relaxed original problem is achieved and proceed to phase two. Otherwise, generate Benders cut and go to step 2). Set $\nu \leftarrow \nu +1$. 
	\item[5):] \textbf{Phase two master problem solution}: Solve the master problem while keeping all the Benders cut generated from phase one.
	\item[6):] \textbf{Unrelaxed subproblem solution}: Enforce the binary constraint for the flow direction $y_{kct}$. Solve the unrelaxed subproblem into optimality and output the optimal solution $\hat{y}_{kct}$. 
	\item[7):] \textbf{Sensitivities generation}:  Fix $y_{kct}=\hat{y}_{kct}$. Solve the subproblem and obtain the dual variables associated with constraints (\ref{fix_P}) and (\ref{fix_de}).
	\item[8):] \textbf{Convergence check}: If $|Z_{up}^{(\nu)}-Z_{down}^{(\nu)}|/|Z_{down}^{(\nu)}|\le \epsilon$, the optimal solution is obtained. Otherwise, generate Benders cut and go to step 5). Set $\nu \leftarrow \nu +1$. 
\end{itemize}

The two phase approach is an efficient method to accelerate convergence of Benders Decomposition \cite{mybibb:bd_twophase}. In phase one, we solve the master problem with the relaxed subproblem at optimality. All the Benders cuts generated in phase one are valid for the original problem. The reason is that the relaxed subproblem provides a lower bound on the original subproblem so that it will also generate a valid lower bound for $\alpha$ \cite{mybibb:valid_cut}. In addition, the objective value obtained from phase one provides a lower bound for the original problem, which can be used to evaluate the quality of the final solution. In phase two, the generation of the Benders cut is heuristic because it involves fixing the binary variable $y_{kct}$. Although it cannot ensure a global optimum, our case studies show that the solution obtained is very close to the lower bound of the original problem. Therefore, the solution is of high quality from an engineering point of view.  

\begin{figure*}[!htb]
	%\captionsetup{font=footnotesize}
	\centering
	\includegraphics[width=0.8\textwidth]{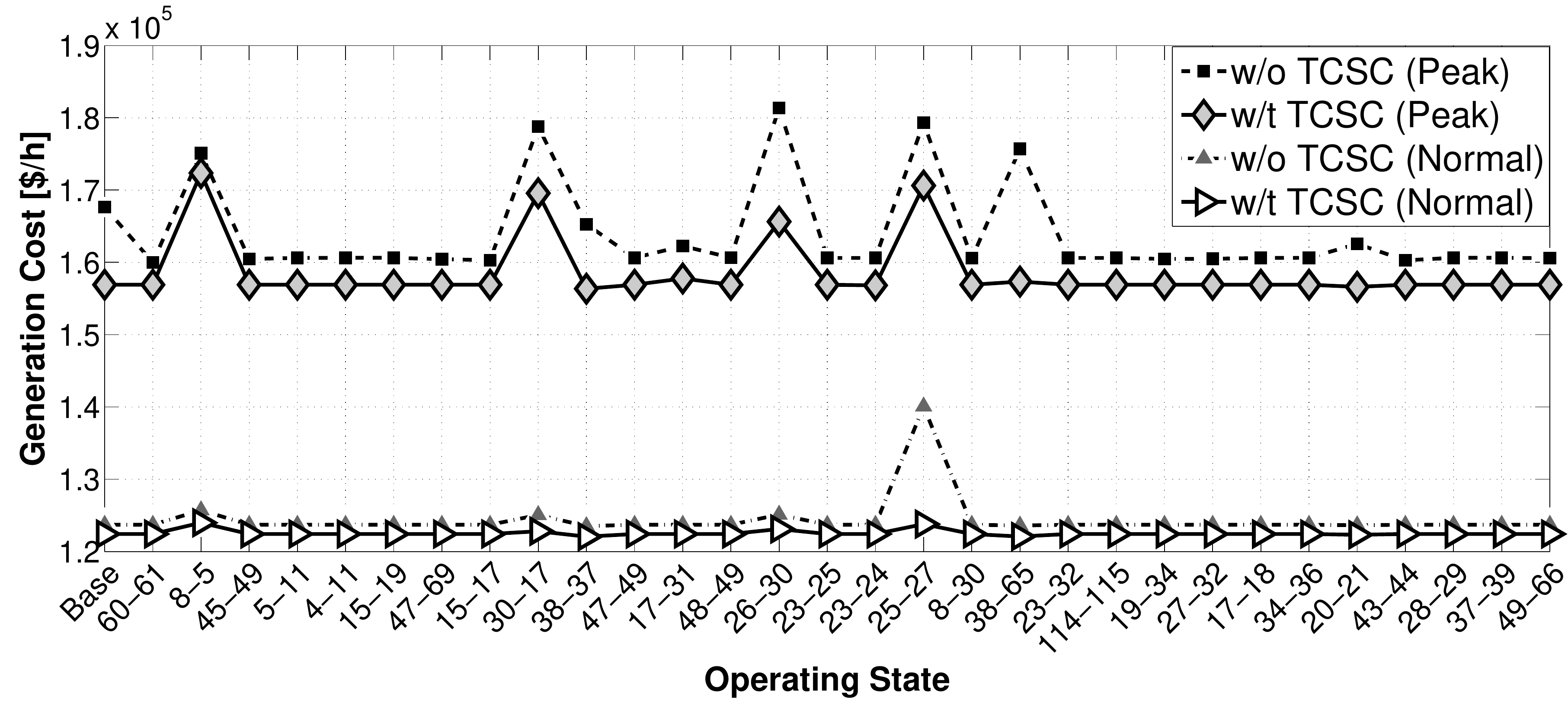}
	\caption{Hourly generation cost for peak and normal load level.}
	\label{hourly_cost}
\end{figure*}     

\section{Numerical Case Studies}
\label{case_study}
The IEEE 118-bus and the Polish 2383-bus system are selected to test the effectiveness of our planning model and the solution approach. The system data is provided in the MATPOWER package \cite{mybibb:MATPOW}. There is only one load pattern defined for these standard systems. For the IEEE 118-bus system, we treat the given loads as the normal load levels. The peak loads are 1.2 times the given loads and the lower loads are 80\% of the given loads. For the Polish system, the provided load data in \cite{mybibb:MATPOW} is the winter peak so we treat the normal loads and low loads as 80\% and 60\% of the given loads, respectively. All simulations are performed on a personal laptop with an Inter Core(TM) i5-2400M CPU @ 2.30 GHz and 4.00 GB of RAM. The problem is modeled by using the MATLAB toolbox YALMIP \cite{mybibb:YALMIP} with CPLEX \cite{mybibb:CPLEX} as the solver.

In this study, we investigate the allocation strategy for one typical VSR: TCSC. The allowable compensation range of TCSC varies from -70\% to +20\% of the corresponding line reactance \cite{mybibb:TCSC_cost_recovery}. Thus, the physical limits for $b^V_{k}$ are $-\frac{1}{6x_k} \le b^V_k \le \frac{7}{3x_k}$. The value of $M_k$ in (\ref{if3}) and (\ref{if4}) is selected as $|\frac{7}{3x_k}\theta_k^{\max}|$. The investment cost of the TCSC depends on its lifespan and capacity rating. The annual investment cost $A_I$ is converted from the total investment cost by using the interest rate and life span of the TCSC \cite{mybibb:GA_FACT_market,mybibb:PSO_SQP_FACTS}. The interest rate is selected to be 5\% and the life span of the TCSC is 5 years \cite{mybibb:TCSC_cost_recovery}.

\subsection{IEEE 118-Bus System}
The IEEE-118 bus system has 118 buses, 19 generators, 177 transmission lines and 9 transformers. The total load at the peak level is 4930 MW and the generation capacity is 6466 MW. The thermal flow limits are decreased artificially to create congestion. In practice, it is unnecessary to consider every transmission line as candidate location to install FACTS device due to some physical or economic limitations. Hence, we first perform a preliminary simulation based on the sensitivity approach proposed in \cite{mybibb:TCSCsens} to obtain 30 TCSC candidate locations. In addition, based on the congestion severity \cite{mybibb:automatic_contingency}, we consider 30 contingencies in the planning model so the number of operating states in this test case is 93. 

Table \ref{comparison} provides a comparison of the non-decomposed approach, the MBD approach and proposed Benders algorithm. The non-decomposed model indicates solving the complete model in Section \ref{Optimization_model} directly \cite{mybibb:investment_naps}. For such a large optimization problem, it may take excessive time to get a solution within the default mipgap (0.01\%) in CPLEX. For the sake of comparison, we simply seek a solution within a given computation time. As shown in the table, the total planning cost for the non-decomposed model is \$1099.59 M with an mipgap 1.47\% after 3 hours. In addition, two TCSCs are selected to be installed in the system. The results for the MBD approach show that five TCSCs should be installed in the system and the total planning cost is \$1090.03 M. The computation time decreases significantly using MBD, requiring only 315.86 s. The proposed Benders algorithm suggest to install six TCSCs in the system and the total planning cost is \$1,088.21M. The lower bound from phase one for this test system is \$1,087.20M, indicating that the solution obtained by the proposed approach is close enough to the global optimal solution. \textcolor{black}{The required computational time is 244.81 s. As compared with MBD, the total planning cost is decreased by 0.20\% and the computational time is reduced by 22.22\%}  

\begin{table}[!htb]
	\centering
	\processtable{Comparison of the Investment Results for IEEE 118-Bus System \label{comparison}}
	%\caption{Comparison of the Investment Results for IEEE 118-Bus System}
	{\tabcolsep7pt
		\begin{tabular}{cccc}
			\hline
			\multirow{2}{*}{Approach}&Non-decomposed &\multirow{2}{*}{MBD \cite{mybibb:bingqian_hu}}& \multirow{2}{*}{Proposed BD}    \\
			&\cite{mybibb:investment_naps}&&  \\
			\hline
			TCSC&\multirow{3}{*}{(26-30),(30-38)}&(17-31),(20-21)& (17-31),(20-21)   \\
			Locations&&(26-30),(22-23)&(21-22),(26-30)    \\
			($i-j$)&&(30-38)&(22-23),(30-38) \\
			\hline
			Investment&\multirow{2}{*}{1.64}&\multirow{2}{*}{2.68}&\multirow{2}{*}{2.94}   \\
			{[million \$]}&&&            \\
			\hline
			Total Cost  &\multirow{2}{*}{1099.59}&\multirow{2}{*}{1090.03}&\multirow{2}{*}{1088.21}   \\
			{[million \$]}&&&  \\
			\hline
			CPU time &3.00 [hours]&315.86 [s]&244.81 [s]    \\
			\hline
	\end{tabular}}{}
\end{table}  

Fig. \ref{hourly_cost} shows hourly generation cost for each operating state under the peak and normal load level. The generation cost reduction can be observed for all the operating states by installing TCSCs into the network. The hourly generation cost for the base case during the peak load level is \$167,653 per hour without any TCSC. This cost decreases to \$156,907 per hour with the installation of six TCSCs. The cost reduction is mainly due to the congestion relief which enables more power to be delivered from cheaper generators. It can also be seen that the generation cost reduction at the normal load level is not as much as that in the peak load level for all the operating states except for contingency (25-27).

\begin{figure}[!htb]
	%\captionsetup{font=footnotesize}
	\centering
	\includegraphics[width=0.65\textwidth]{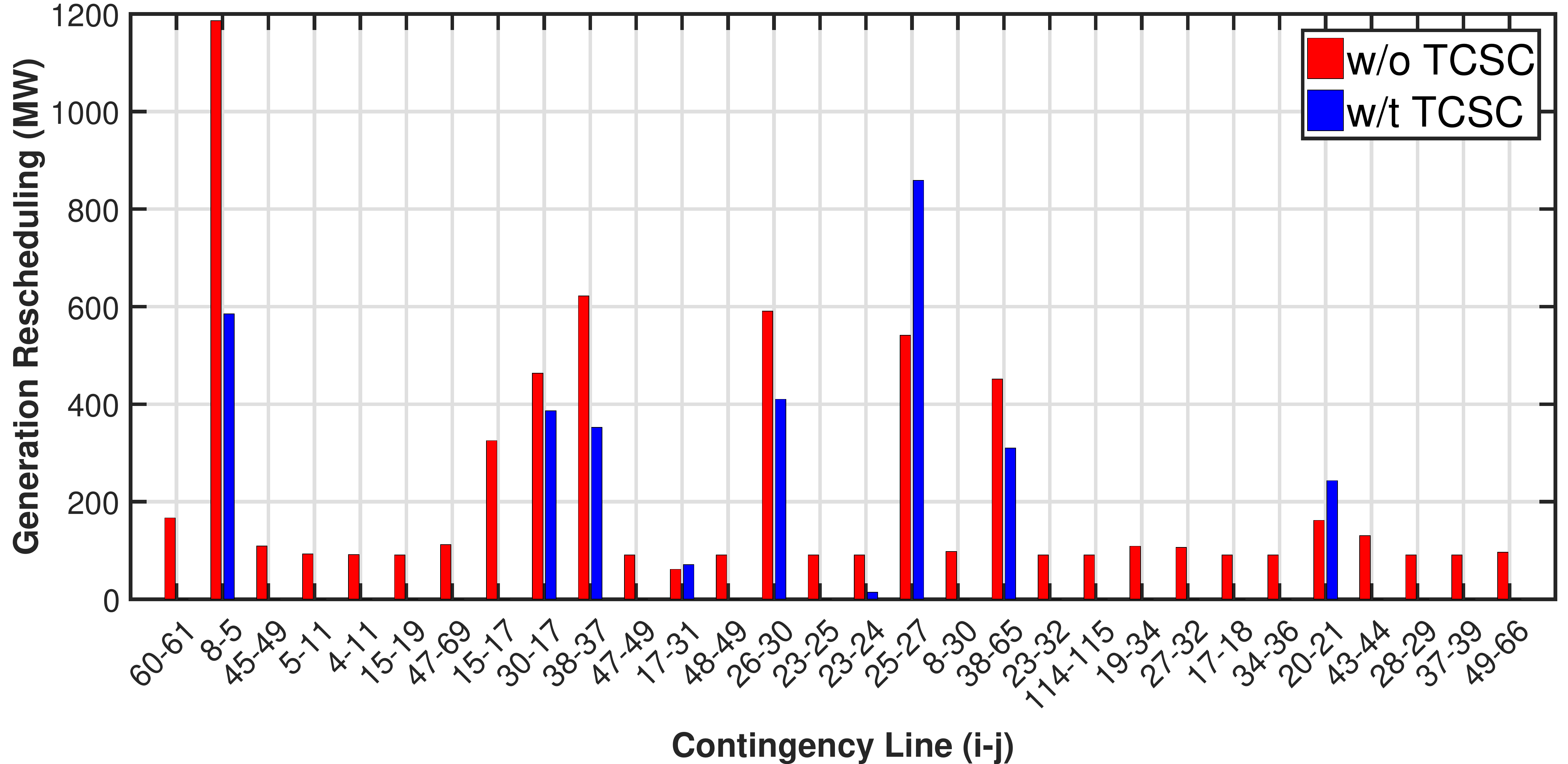}
	\caption{Generation rescheduling under different contingencies for the peak load level.}
	\label{resch}
\end{figure}

\begin{figure}[!htb]
	%\captionsetup{font=footnotesize}
	\centering
	\includegraphics[width=0.65\textwidth]{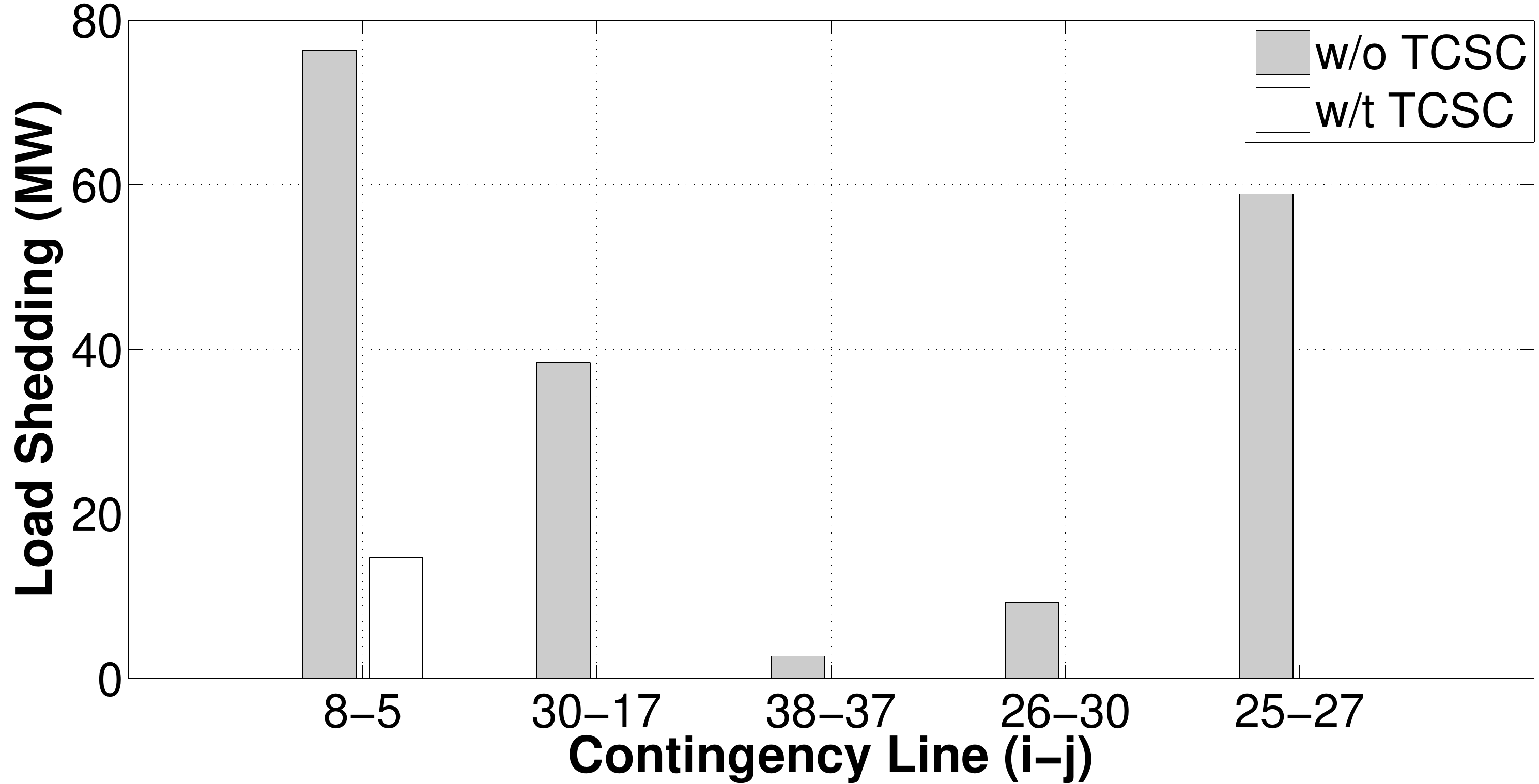}
	\caption{Load shedding under different contingencies for peak load level.}
	\label{LS}
\end{figure}

Fig. \ref{resch} shows generation rescheduling under different contingencies for the peak load level. Fig. \ref{LS} gives the load shedding at the peak load level for the five contingencies which involve load shedding. Note that there is no involuntary load shedding for the operating states under normal and low load level. From Fig. \ref{resch}, it can be seen that the amount of generation rescheduling decreases for the majority of operating states. The largest reduction occurs under contingency (8-5) where the amount of generation rescheduling decreases from 1200 MW to about 600 MW. Under contingency (25-27), the rescheduling amount increases by about 300 MW with TCSC. However, about 60 MW load shedding can be avoided in that contingency as shown in Fig. \ref{LS}. This indicates that the installation of TCSC enable cheaper ways, such as, rescheduling, to reduce load shedding. As can be seen in Fig. \ref{LS}, the load shedding for contingency (30-17), (38-37), (26-30) and (25-27) are eliminated with TCSC. For the most severe contingency (8-5), the load shedding decreases from 76.36 MW to 15.30 MW.   

Table \ref{annual_saving_118} compares the annual cost for the case with and without TCSC. We categorize the cost into four groups: 1) generation cost; 2) generation rescheduling cost; 3) load shedding cost; and 4) investment cost in TCSC. Except for the investment in TCSC, it can be seen that cost decreases in all the other categories with the installation of TCSCs. The annual reduction for the total planning cost is about \$36.58 M, which approximately accounts for 3.25\% of the annual planning cost.

\begin{table}[!htb]
	\centering
	\processtable{Annual Cost with and without TCSC for IEEE 118-Bus System \label{annual_saving_118}}
	{\tabcolsep8pt
		\begin{tabular}{ccc}
			\hline
			\multirow{2}{*}{Cost Category}&\multicolumn{2}{c}{Annual Cost [million \$]}     \\
			%\cline{2-3}
			&w/o TCSC&w/t TCSC      \\
			\hline
			Generation cost in normal state&1048.31&1018.74   \\
			\hline
			Generation cost in contingency&66.58&65.30     \\
			\hline
			Rescheduling cost&1.09&0.54           \\
			\hline
			Load shedding cost&8.81&0.70      \\
			\hline
			Investment in TCSC&-&2.94     \\
			\hline
			Total cost&1124.79&1088.21      \\
			\hline
	\end{tabular}}{}
\end{table}

\begin{figure}[!htb]
	%	\captionsetup{font=footnotesize}
	\centering
	\includegraphics[width=0.65\textwidth]{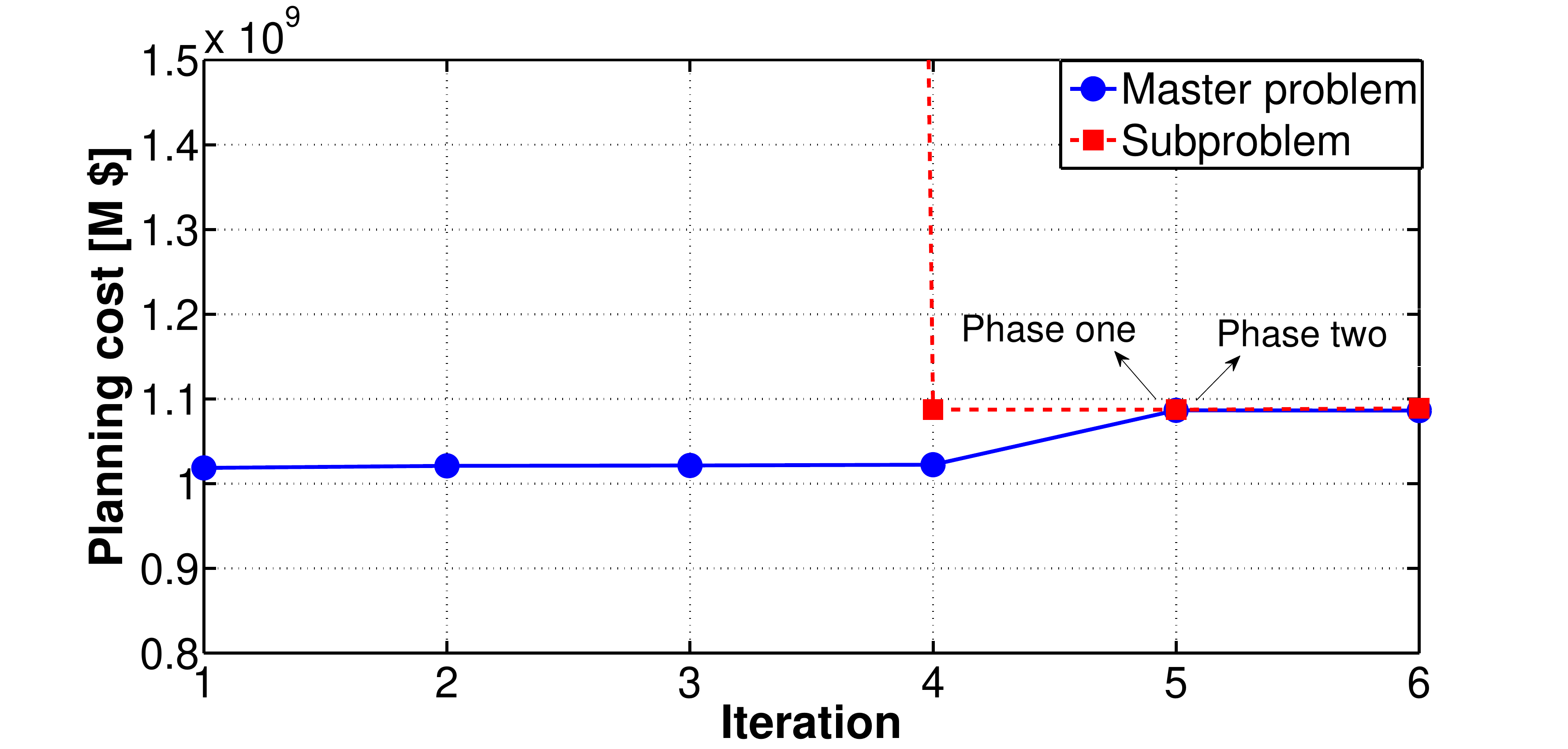}
	\caption{\textcolor{black}{Evolution of the proposed Benders algorithm for IEEE 118-bus system.}}
	\label{iteration_118}
\end{figure}

\textcolor{black}{Fig. \ref{iteration_118} depicts the convergence on the IEEE 118-bus system. In the first three iterations, the penalty terms in (\ref{sub_obj}) are not zero so the objective value of subproblem is very large. We decrease the range of y-axis to improve the readability of the figure. As can be observed from the figure, it takes five iterations for the problem in phase one to converge. After that, only one iteration is needed for the problem in phase two to converge.} 

\subsection{Polish System}
The Polish system includes 2,383 buses, 327 generators, 2,728 transmission lines and 168 transformers. The total load at the peak level is 20,465 MW and the generation capacity is 29,594 MW. Based on sensitivity method, we select 50 candidate locations to install TCSC. We consider 60 contingencies so the number of operating states for this test system is 183.  

\textcolor{black}{Table \ref{comparison_2383} gives a comparison of the MBD approach and the proposed algorithm. The non-decomposed approach \cite{mybibb:investment_naps} is not included in the table since the CPLEX solver gives no solution after five hours on this system. The convergence tolerance $\epsilon$ is selected to be 0.35\%. In addition, a time limit of three hours is set for each algorithm. At the end of each iteration, if the summation of the computational time for master problem and sub-problem exceeds three hours, the algorithm will terminate. As can be seen from the table, the MBD suggests to install 13 TCSCs and the total planning cost is \$10796.54 M. Moreover, it takes 3.13 hours for the MBD to achieve a gap of 1.12\%. Fifteen TCSCs are selected to be installed in the system by using the proposed algorithm. The planning cost is \$10778.41 M and the computational time is decreased to 1.5 hours. The proposed algorithm demonstrates better performance in relieving the computational burden when applied to practical large scale system as compared with MBD. }
\begin{table}[!htb]
	\centering
	\processtable{\textcolor{black}{Comparison of the Investment Results for Polish System} \label{comparison_2383}}
	%\caption{Comparison of the Investment Results for IEEE 118-Bus System}
	{\tabcolsep7pt
		\begin{tabular}{ccc}
			\hline
			\multirow{2}{*}{Approach}&\multirow{2}{*}{MBD \cite{mybibb:bingqian_hu}}& \multirow{2}{*}{Proposed BD}    \\
			&&  \\
			\hline
			&(1507-1303),(123-111)& (29-13),(7-4),(10-3)   \\
			&(432-356),(920-821)&(1342-1301),(1948-1649)    \\
			TCSC&(7-4),(1181-1190)&(432-356),(920-821) \\
			Locations&(1467-834),(10-5)&(395-334),(493-306)  \\
			($i-j$)&(612-413),(114-131)&(11-4),(152-66)  \\
			&(578-366),(1489-1431)&(612-413),(1489-1431)  \\
			&(29-7)&(833-1230),(1055-1079)  \\
			\hline
			Investment&\multirow{2}{*}{12.34}&\multirow{2}{*}{12.43}   \\
			{[million \$]}&&            \\
			\hline
			Total Cost  &\multirow{2}{*}{10796.54}&\multirow{2}{*}{10778.41}   \\
			{[million \$]}&&  \\
			\hline
			CPU time &3.13 [h]&1.50 [h]    \\
			\hline
	\end{tabular}}{}
\end{table}  

Table \ref{annual_saving_2383} shows the comparison of the planning cost for the case with and without TCSC. The annual savings for the Polish system is about \$64.5 M.

\begin{table}[!htb]
	\centering
	\processtable{Annual Cost with and without TCSC for the Polish 2383-Bus System \label{annual_saving_2383}}
	{\tabcolsep8pt
		\begin{tabular}{c c c}
			\hline
			\multirow{2}{*}{Cost Category}&\multicolumn{2}{c}{Annual Cost [million \$]}     \\
			%\cline{2-3}
			%\hline
			&w/o TCSC&w/t TCSC      \\
			\hline
			Generation cost in normal state&9527.18&9464.42   \\
			\hline
			Generation cost in contingency&1299.87&1291.18     \\
			\hline
			Rescheduling cost&6.28&4.83           \\
			\hline
			Load shedding cost&9.59&5.55      \\
			\hline
			Investment in TCSC&-&12.43     \\
			\hline
			Total cost&10842.93&10778.41      \\
			\hline
	\end{tabular}}{}
\end{table}

\begin{figure}[!htb]
	%	\captionsetup{font=footnotesize}
	\centering
	\includegraphics[width=0.65\textwidth]{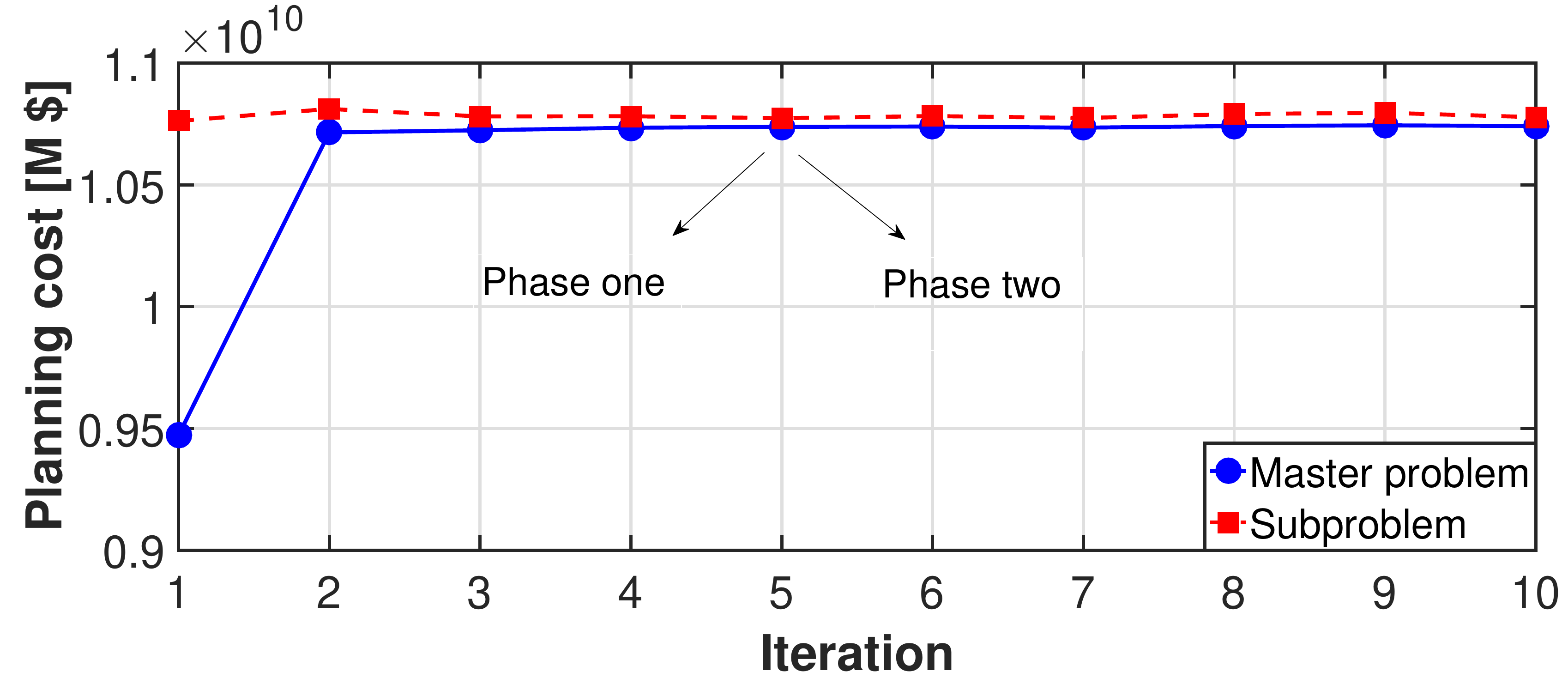}
	\caption{Evolution of the proposed Benders algorithm for Polish system.}
	\label{iteration}
\end{figure}

Fig. \ref{iteration} illustrates the iteration process of the proposed Benders algorithm. It can be seen that after 5 iterations, the problem in phase one converges to the lower bound of  \$10,774.27 M. Then it takes another 5 iterations for the problem in phase two to converge. 

\subsection{Discussion}
\textcolor{black}{The proposed approach shows better performance than the MBD in determining the optimal locations of TCSC. Although the difference of the total planning cost between these two approaches is not significant, the computational time is greatly reduced by using the proposed approach, especially for the practical large scale system. Another advantage of the proposed approach is that it can provide a lower bound for the complete planning model, i.e., the objective value of phase one. This information can be leveraged to evaluate the quality of the obtained solution.} 

%\textcolor{blue}{Note that the proposed planning model adopts DC power flow, which ignores the power loss and reactive power. Thus, the proposed model is suitable for %preparatory power network design. The obtained results need to be further investigated by using AC power flow model. After that, the power losses, different security and %stability index should be analyzed and evaluated.}

\textcolor{black}{The focus of this paper is to determine the optimal locations of series FACTS devices based on cost. The obtained solution provides a good trade-off between FACTS investment cost and system operation cost. Note that the proposed planning model adopts a DC power flow, which ignores power loss and reactive power. Thus, the proposed model is suitable for preparatory power network design. For more detailed security analysis, the obtained results can be analyzed by using an AC power flow model. In that analysis, the power losses and various security indices \cite{mybibb:security_index} can be computed and analyzed.}

\section{Conclusion}
\label{conclusions}
This paper proposes a planning model to allocate VSR considering different operating conditions and critical $N-1$ contingencies. The original planning model is a large scale MINLP model. A reformulation is introduced to transform the MINLP model into a MILP model. To further reduce the computation burden, a two phase Benders decomposition is proposed. The solution obtained is not guaranteed to be a global optimum but analysis indicates the solution is near optimal. Case studies on the IEEE 118-bus and the Polish system demonstrate the performance of the proposed algorithm. The simulation results show that the generation cost for both the normal operating states and contingency states can be reduced with the installation of VSR. In addition, the cost reductions can be observed in the generation rescheduling and involuntary load shedding following contingencies.    

\section{Acknowledgments}
This work was supported in part by ARPAe (Advanced Research Projects Agency Energy), and in part by the Engineering Research Center Program of the National Science Foundation and the Department of Energy under NSF Award Number EEC-1041877 and the CURENT Industry Partnership Program.

\end{document}